\documentclass[10pt]{article}
\usepackage{amsfonts}
\usepackage{amssymb}
\usepackage{amsthm}
\usepackage[centertags]{amsmath}
\usepackage{newlfont}
\usepackage{graphicx}

\newfont{\bb}{msbm10 at 10pt}
\def\r{\hbox{\bb R}}
\def\h{\hbox{\bb H}}
\def\s{\hbox{\bb S}}
\newcommand{\T}{\hbox{\bf t}}
\newcommand{\N}{\hbox{\bf n}}
\newcommand{\B}{\hbox{\bf b}}
\newcommand{\C}{\hbox{\bf c}}
\newtheorem{theorem}{Theorem}[section]

\newtheorem{proposition}[theorem]{Proposition}

\newtheorem{lemma}[theorem]{Lemma}
\newtheorem{corollary}[theorem]{Corollary}

\begin{document}

\title{Linear Weingarten surfaces in Euclidean and hyperbolic space}
\author{
 Rafael L\'opez\footnote{Partially
supported by MEC-FEDER
 grant no. MTM2007-61775 and
Junta de Andaluc\'{\i}a grant no. P06-FQM-01642.}}
\date{}

\maketitle
\begin{center}\emph{This paper is dedicated to Manfredo do Carmo in admiration for his mathematical achievements and his  influence on the field of differential geometry of surfaces}
\end{center}
\begin{abstract} In this paper we review some author's results about Weingarten surfaces in Euclidean space $\r^3$ and hyperbolic space $\h^3$. We stress here in the search of examples of linear Weingarten surfaces that satisfy a certain geometric property. First, we consider Weingarten surfaces in $\r^3$ that are foliated by circles, proving that the surface is rotational, a Riemann example or a generalized cone. Next we classify rotational surfaces in $\r^3$ of hyperbolic type showing that there exist surfaces that are complete. Finally, we study linear Weingarten surfaces in $\h^3$ that are invariant by a group of parabolic isometries, obtaining its classification.
\end{abstract}

\emph{MSC:}  53C40, 53C50

\emph{Keywords}:  Weingarten surface; cyclic surface; Riemann examples; parabolic surface.

%%%%%%%%%%%%%%%%%%%%%%%%%%%%%%%%%%%%%%%%%%%%%%%
\section{Statement of results}
%%%%%%%%%%%%%%%%%%%%%%%%%%%%%%%%%%%%%%%%%%%%%%%

A surface $S$ in Euclidean space $\r^3$ or hyperbolic space $\h^3$ is called a {\it
Weingarten surface} if there is some smooth relation $W(\kappa_1,\kappa_2)=0$ between its
two principal curvatures $\kappa_1$ and $\kappa_2$. In particular, if $K$ and $H$  denote respectively  the Gauss curvature and the mean curvature of  $S$, $W(\kappa_1,\kappa_2)=0$ implies a relation $U(K,H)=0$. The  classification of  Weingarten surfaces in
the general case  is almost completely open today. After earlier works in the fifties due to Chern, Hopf, Voss, Hartman, Winter, amongst others,   there has been recently a progress in this theory, specially when the Weingarten relation is of type $H=f(H^2-K)$ and $f$ elliptic. In such case, the surfaces satisfy a maximum principle that allows a best knowledge of the shape of such surfaces. These achievements can see, for example,  in \cite{bs,co,gmm,rs,st1,st2}.

The simplest case of functions $W$ or $U$ is that they are linear, that is,
\begin{equation}\label{rlopez:w1}
a\kappa_1+b\kappa_2=c\hspace*{1cm}\mbox{or}\hspace*{1cm}a H+b K=c,
\end{equation}
where $a,b$ and $c$ are constant. Such surfaces are called \emph{linear Weingarten} surfaces. Typical examples of linear Weingarten surfaces are  umbilical surfaces,  surfaces with constant Gauss curvature and surfaces with constant mean curvature.

A first purpose of the present work is to provide  examples of linear Weingarten surfaces that satisfy a certain geometric condition. A first attempt is to consider that the surface is rotational, that is, invariant by a group of isometries that leave fixed-pointwise a geodesic of the ambient space. In such case, equations (\ref{rlopez:w1}) lead to  an ordinary  differential equations and the study is then reduced to finding the profile curve that defines the surface.

A more general family of rotational surfaces are the cyclic surfaces, which were introduced by Enneper in the XIX century.  A \emph{cyclic surface} in Euclidean space $\r^3$ or $\h^3$ is a surface determined by a smooth uniparametric family of  circles.
Thus, a cyclic surface $S$ is a surface  foliated by  circles meaning
that there is a one-parameter family of planes which meet $S$  in these
 circles.  The planes are not assumed parallel, and
if two circles should lie in planes that happen to be parallel, the
circles are not assumed coaxial. Rotational surfaces are examples of cyclic surfaces.

Our first result is motivated by what happens for cyclic  surfaces with constant mean curvature $H$. Recall that the catenoid is the only minimal ($H=0$) rotational surface in $\r^3$. If the surface is not rotational, then the only  cyclic minimal surfaces are a family of examples of periodic minimal surfaces discovered by Riemann, usually called in the literature as Riemann examples \cite{ri}. If the mean curvature $H$ is a non-zero constant, then the only  cyclic surfaces are the surfaces of revolution (Delaunay surfaces) \cite{ni}. In order to find new examples of linear Weingarten surfaces, we pose the following question: do exist non-rotational cyclic surfaces that are linear Weingarten surfaces?

In Section \ref{rlopez:sec2} we prove the following result:

\begin{theorem}\label{rlopez:t1} Let $S$ be a cyclic  surface in
Euclidean space $\r^3$.
\begin{enumerate}
\item If $S$ satisfies a relation of type $a\kappa_1+b\kappa_2=c$, then $S$ is  a surface of revolution or it is  a Riemann  example ($H=0$).
 \item If $S$ satisfies a relation of type $a H+b K=c$, then $S$ is  a surface of revolution or it is  a Riemann  example ($H=0$) or it is  a generalized cone ($K=0$).
     \end{enumerate}
\end{theorem}
Recall that a generalized cone is
a cyclic surface formed by a uniparametric family of $u$-circles whose  centres lie in a
straight line and the radius function is linear on $u$. These surfaces have $K\equiv 0$ and they are the only  non-rotational cyclic surfaces in $\r^3$ with constant Gaussian curvature \cite{lo}.

After Theorem \ref{rlopez:t1}, we focus on Weingarten surfaces of revolution in $\r^3$. The classification of linear Weingarten surfaces strongly depends on
  the sign of $\Delta:=a^2+4bc$. If $\Delta>0$, the surface is said elliptic and satisfies good properties, as for example, a maximum principle: see \cite{gmm,rs}. If $\Delta=0$, the surface is a tube, that is, a cyclic surface where the circles have the same radius. Finally, if $\Delta<0$, the surface is said \emph{hyperbolic} (see \cite{ae}). In Section \ref{rlopez:sec3} we study  hyperbolic rotational surfaces in $\r^3$.  We do an explicit description of the hyperbolic rotational  linear Weingarten. Examples of hyperbolic Weingarten surfaces are the surfaces with constant negative Gaussian curvature $K$: we take $a=0$, $b=1$ and $c<0$ in the right relation (\ref{rlopez:w1}). In contrast to the  Hilbert's theorem  that asserts that do not exist  complete surfaces with constant negative Gaussian
curvature immersed in $\r^3$, we obtain (see Theorem \ref{rlopez:t42}):

\begin{theorem}\label{rlopez:t2} There exists a family  of  hyperbolic linear
Weingarten complete rotational surfaces in $\r^3$ that are non-embedded and periodic.
\end{theorem}

Finally, we are interested in  linear Weingarten surfaces of revolution in hyperbolic space $\h^3$. In hyperbolic space
there exist three types of rotational surfaces. We will study one of them, called  \emph{parabolic surfaces}, that is, surfaces invariant by a group of parabolic isometries of the ambient space.  This was began by do Carmo and Dajczer in the study of rotational surfaces in $\h^3$ with constant curvature \cite{cd} and  works of Gomes, Leite, Mori et al. We will consider problems such as existence, symmetry and behaviour at infinity. As a consequence  of our work, we obtain the following

\begin{theorem} \label{rlopez:t3}  There exist  parabolic complete surfaces in $\h^3$  that satisfy the relation $aH+bK=c$.
\end{theorem}

Part of the results of this work have recently appeared in a series of author's papers:  \cite{lo1,lo2,lo3,lo4}.

%%%%%%%%%%%%%%%%%%%%%%%%%%%%%%%%%%%%%%%%%%%%%%%%%%%%%%%%%%%%%%%%%%%%%%%%%

\section{Cyclic  Weingarten surfaces in $\r^3$}\label{rlopez:sec2}

%%%%%%%%%%%%%%%%%%%%%%%%%%%%%%%%%%%%%%%%%%%%%%%%%%%%%%%%%%%%%%%%%%%%%%%%%

In this section we will study linear Weingarten surfaces that are foliated by a  uniparametric family of circles (cyclic surfaces). In order to show the techniques to get Theorem \ref{rlopez:t1}, we only consider  Weingarten surfaces that satisfy the linear relation $aH+bK=c$. The proof consists into two steps. First, we prove

\begin{theorem}\label{rlopez:t21} Let  $S$ be a  surface that satisfies $aH+bK=c$ and it is foliated by circles
lying in a one-parameter family of planes. Then either $S$ is a
subset of a round sphere or the planes of the foliation  are parallel.
\end{theorem}

\begin{proof} Consider $P(u)$ the set of planes of the foliation, that is, $S=\bigcup_{u\in I}P(u)\cap S$, $u\in I\subset\r$, and such that
$P(u)\cap S$ is a circle for each $u$. Assume that the planes $P(u)$ are not parallel . Then we are going to show that the surfaces is a sphere. The proof follows the same ideas for the case of the constancy of the mean curvature \cite{ni}. Let $\Gamma$ be an orthogonal curve to the foliation planes, that is, $\Gamma'(u)\bot P(u)$.
 If $\{\T,\N,\B\}$ denotes the usual Frenet trihedron of $\Gamma$, the surface $S$ is locally parametrized by
\begin{equation}\label{rlopez:para}
X(u,v)=\C(u)+r(u)(\cos{v}\ \N (u)+\sin{v}\ \B (u)),
\end{equation}
where $r=r(u)>0$ and $\C=\C(u)$ denote respectively
the radius and centre of each circle $P(u)\cap S$. We compute the mean curvature and the Gauss curvature of $S$ using the usual local formulae
$$H=\frac{eG-2f F+gE}{2(EG-F^2)},\hspace*{1cm}
K=\frac{eg-f^2}{EG-F^2}.$$
Here $\{E,F,G\}$ and $\{e,f,g\}$ represent  the coefficients of the first and second fundamental form, respectively.
 Then the relation $aH+bK=c$ writes in terms of the curve $\Gamma$. Using the Frenet equations of $\Gamma$, we are able to
 express  the relation $aH+bK=c$ as a trigonometric polynomial on $\cos{(n v)}$ and $\sin{(nv)}$:
$$A_0+\sum_{n=1}^8 \bigg(A_n(u) \cos{(nv)} +B_n(u)\sin{(n v)}\bigg)=
0,\hspace*{.5cm}u\in I, v\in [0,2\pi].$$
Here $A_n$ and $B_n$ are smooth functions on $u$. Because the functions  $\cos{ (n v)}$ and $\sin{(nv)}$ are independent, all coefficient functions $A_n,B_n$ must be zero. This leads to a set of equations, which we wish to solve. Because the curve $\Gamma$ is not a straight line, its curvature $\kappa$ does not vanish.

The proof consists into the explicit computation of the coefficients $A_n$ and $B_n$ and solving $A_n=B_n=0$. The proof program begins with the equations $A_8=0$ and $B_8=0$, which yields relations between the geometric quantities of the curve $\Gamma$. By using these data, we follow with equations $A_7=B_7=0$ and so on, until to arrive with $n=0$. The author
 was able to obtain the results using the  symbolic  program
Mathematica to check his work: the computer was used in each calculation several times, giving understandable expressions of the coefficients $A_n$ and $B_n$. Finally, we achieve to show that $X$ is a parametrization of a round sphere.

\end{proof}

Once proved Theorem \ref{rlopez:t21}, the following step consists to conclude that either the
circles of the foliation must be coaxial (and the surface is rotational) or that $K\equiv 0$ or $H\equiv 0$. In the latter cases, the
 Weingarten relation (\ref{rlopez:w1}) is trivial in the sense that $a=c=0$ or $b=c=0$.

\begin{theorem} \label{rlopez:t22} Let  $S$ be a  cyclic surface that satisfies $aH+bK=c$. If the foliation planes are parallel,
then  either $S$ is a surface of revolution or $a=c=0$ or $b=c=0$.
\end{theorem}

\begin{proof} After an isometry of the ambient space, we parametrize $S$ as
$$X(u,v)=(f(u),g(u),u)+r(u)(\cos{v},\sin{v},0),$$
where  $f, g$ and $r$ are smooth functions on  $u$, $u\in I\subset\r$  and $r(u)>0$ denotes the
radius of each circle of the foliation. Then   $S$ is a surface of
revolution if and only if $f$ y $g$ are constant functions. Proceeding similarly as in the proof of  Theorem \ref{rlopez:t21}, Equation $aH+bK=c$ is equivalent to an expression
$$\sum_{n=0}^{8} \bigg(A_n(u) \cos{(n v)}+B_n(u)\sin{(n v)}\bigg)=0.$$
Again, the functions $A_n$ and $B_n$  must vanish on $I$. Assuming that the surface is not rotational, that is, $f'(u)g'(u)\not=0$ at some point $u$, we conclude that $a=c=0$ or $b=c=0$.
\end{proof}

Recall what happens in the latter cases. The computation of $H\equiv 0$ and $K\equiv 0$ gives
\begin{equation}\label{rlopez:riemann}
f''=\lambda r^2,\hspace*{.5cm}g''=\mu r^2,\hspace*{.5cm}
1+(\lambda^2+\mu^2) r^4+r'^2-r r''=0,
\end{equation}
and
\begin{equation}\label{rlopez:cone}
f''=g''=r''=0,
\end{equation}
respectively. If (\ref{rlopez:riemann}) holds, we have the equations that describe the  Riemann  examples ($\lambda^2+\mu^2\not=0$) and
the catenoid ($\lambda=\mu=0$). In the case (\ref{rlopez:cone}), the surface $S$ is a generalized cone.

As a consequence of the above Theorems \ref{rlopez:t21} and \ref{rlopez:t22}, we obtain  Theorem \ref{rlopez:t1} announced in the introduction of this work. Finally, the previous results allow us to give a characterization of Riemann examples and generalized cones in the class of linear Weingarten surfaces.

 \begin{corollary} Riemann examples   and generalized cones  are the only non-rotational cyclic surfaces that satisfy a  Weingarten relation of type $aH+bK=c$.
\end{corollary}

%%%%%%%%%%%%%%%%%%%%%%%%%%%%%%%%%%%%%%%%%
\section{Hyperbolic linear Weingarten surfaces in $\r^3$}\label{rlopez:sec3}
%%%%%%%%%%%%%%%%%%%%%%%%%%%%%%%%%%%%%%%%%%%%%%%%%%

We consider  surfaces $S$ in Euclidean space that satisfy the relation
\begin{equation}\label{rlopez:w31}
a\ H+b\ K=c
\end{equation}
where $a$, $b$ and $c$ are constants under the relation $a^2+4bc<0$. These  surfaces are called \emph{hyperbolic linear Weingarten surfaces}.
In particular, $c\not=0$, which can be assumed to be $c=1$. Thus the condition $\Delta<0$ writes now as $a^2+4b<0$. In this section, we study these surfaces in the class of surfaces of revolution.  Equation (\ref{rlopez:w31}) leads to an ordinary  differential equation that describes the generating curve $\alpha$ of the surface. Without loss of generality, we assume $S$ is a rotational surface whose axis is the $x$-axis. If  $\alpha(s)=(x(s),0,z(s))$ is arc-length parametrized and the surface is given by $X(s,\phi)=(x(s), z(s)
\cos\phi, z(s)\sin\phi)$, then (\ref{rlopez:w31}) leads to
\begin{equation}\label{rlopez:w32}
a\ \frac{\cos\theta(s)-z(s)\theta'(s)}{2 z(s)}-b\ \frac{\cos\theta(s)\theta'(s)}{z(s)}=1,
\end{equation}
where  $\theta=\theta(s)$ the angle function that makes the
velocity $\alpha'(s)$ at $s$  with the $x$-axis, that is,
$\alpha'(s)=(\cos\theta(s),0,\sin\theta(s))$. The curvature of the
planar curve $\alpha$ is given by $\theta'$. In this section, we
 discard the trivial cases in (\ref{rlopez:w31}), that is, $a=0$ (constant Gauss curvature) and $b=0$ (constant mean curvature).

The generating curve $\alpha$ is then described by the solutions of the O.D.E.
\begin{equation}\label{rlopez:eq1}
 \left\{\begin{array}{lll}
 x'(s)&=& \displaystyle \cos\theta(s)\\
 z'(s)&=&\displaystyle \sin\theta(s)\\
 \theta'(s)&=&\displaystyle \frac{a\cos\theta(s)-2z(s)}{a z(s)+2b \cos\theta(s)}
\end{array}
\right.
\end{equation}
Assume  initial conditions
\begin{equation}\label{rlopez:eq2}
x(0)=0,\hspace*{.5cm}z(0)=z_0,\hspace*{.5cm}\theta(0)=0.
\end{equation}
Without loss of generality, we can choose the parameters $a$ and $z_0$ to have the same sign: in our case, we take to be positive numbers.

A first integral of (\ref{rlopez:eq1})-(\ref{rlopez:eq2}) is given by
\begin{equation}\label{rlopez:first}
z(s)^2-a z(s)\cos\theta(s)-b\cos^2\theta(s)-(z_0^2-az_0-b)=0.
\end{equation}
By the uniqueness of solutions, any solution $\alpha(s)=(x(s),0,z(s))$  of (\ref{rlopez:eq1})-(\ref{rlopez:eq2}) is symmetric with respect to the line $x=0$.

In view of (\ref{rlopez:eq2}), the value of $\theta'(s)$ at $s=0$ is $\theta'(0)=\displaystyle \frac{a-2z_0}{a z_0+2b}$. Our study depends on the sign of $\theta'(0)$. We only consider the case
\begin{equation}\label{rlopez:c4}
z_0>\frac{-2b}{a}
\end{equation}
which implies that $z_0>a/2$. The denominator in the third equation of (\ref{rlopez:eq1}) is positive since it does not vanish and at $s=0$, its value is $a z_0+2b>0$. As $z_0>a/2$, the numerator in (\ref{rlopez:eq1}) is negative. Thus we conclude that the function $\theta'(s)$ is negative
anywhere.

From (\ref{rlopez:first}), we write the function $z=z(s)$ as
\begin{equation}\label{rlopez:zeta}
z(s)=\frac12\left(a\cos\theta(s)+\sqrt{(a^2+4b)\cos^2\theta(s)+4(z_0^2-az_0-b)}\right).
\end{equation}

\begin{lemma}\label{rlopez:le1}
 The maximal interval of the solution $(x,z,\theta)$ of  (\ref{rlopez:eq1})-(\ref{rlopez:eq2}) is $\r$.
\end{lemma}

\begin{proof} The result follows if we prove that the derivatives $x', z'$ and $\theta'$ are bounded. In view of (\ref{rlopez:eq1}), it suffices to show it for $\theta'$ (recall that $\theta'(s)<0$). We are going to find a negative number $m$  such that $m\leq
\theta'(s)$ for all $s$. To be exact,  we show the existence of constants $\delta$ and $\eta$
independent on $s$, with
 $\eta<0<\delta$,  such that
\begin{equation}\label{rlopez:delta}
az(s)+2b\cos\theta(s)\geq\delta\hspace*{.5cm}\mbox{and}\hspace*{.5cm}
a\cos\theta(s)-2z(s)\geq\eta.
\end{equation}
Once proved this, it follows from (\ref{rlopez:eq1}) that
\begin{equation}\label{rlopez:m1}
\theta'(s)\geq \frac{\eta}{\delta}:=m.
\end{equation}
Define the
function $f(z_0):=z_0^2-az_0-b$.  The function $f$ is strictly increasing on $z_0$ for
$z_0>a/2$. Using that $a^2+4b<0$, we have
$\frac{a}{2}<\frac{-2b}{a}$. As $z_0$ satisfies (\ref{rlopez:c4}),   there exists  $\epsilon>0$  such that
$$z_0^2-a z_0-b=f(-\frac{2b}{a})+\epsilon=\frac{b(a^2+4b)}{a^2}+\epsilon.$$
From (\ref{rlopez:zeta}),
$$z(s)\geq\frac12\left(a\cos\theta+
\sqrt{(a^2+4b)\cos^2\theta+\frac{4b}{a^2}(a^2+4b)+4\epsilon}\right)\geq
\frac12 \left(a\cos\theta-\frac{a^2+4b}{a}+\epsilon'\right),$$
for a
certain positive number $\epsilon'$. By using that $a^2+4b<0$ again, we have
$$az(s)+2b\cos \theta(s)\geq \frac{a^2+4b}{2}(\cos\theta(s)-1)+\frac{a}{2}\epsilon'\geq \frac{a}{2}\epsilon':=\delta.$$
By using (\ref{rlopez:zeta}) again, we obtain
$$a\cos\theta(s)-2z(s)\geq-\sqrt{(a^2+4b)\cos^2(s)\theta+4f(z_0)}\geq -2\sqrt{f(z_0)}:=\eta,$$
which concludes the proof of this lemma.
\end{proof}

\begin{lemma} \label{rlopez:le2} For each solution $(x,z,\theta)$ of  (\ref{rlopez:eq1})-(\ref{rlopez:eq2}), there exists $M<0$ such that
$\theta'(s)<M$.
\end{lemma}
\begin{proof}
It suffices if we prove that there
exist $\delta_2,\eta_2$, with $\eta_2<0<\delta_2$ such that
$$az(s)+2b\cos\theta(s)\leq\delta_2\hspace*{.5cm}\mbox{and}\hspace*{.5cm}
a\cos\theta(s)-2z(s)\leq\eta_2,$$
since (\ref{rlopez:eq1}) yields $\theta'(s)\leq \delta_2/\eta_2:=M$.  Using (\ref{rlopez:zeta}), we have
\begin{eqnarray*}
az(s)+2b\cos\theta(s)&=&\frac12\left((a^2+4b)\cos\theta(s)+a\sqrt{(a^2+4b)\cos^2\theta(s)+4f(z_0)}\right)\\
&\leq& a\sqrt{f(z_0)}:=\delta_2.
\end{eqnarray*}
On the other hand,
\begin{eqnarray*}
a\cos\theta(s) -2z(s)&=&-\sqrt{(a^2+4b)\cos^2\theta(s)+4f(z_0)}\\
&\leq &-\sqrt{(a^2+4b)+4f(-2b/a)}:=\eta_2.
\end{eqnarray*}
\end{proof}
Lemma \ref{rlopez:le2} implies that $\theta(s)$ is strictly decreasing with
$$\lim_{s\rightarrow\infty}\theta(s)=-\infty.$$
Since Lemma \ref{rlopez:le1} asserts that any  solution is defined for any $s$, put  $T>0$ the first number such that $\theta(T)=-2\pi$.
We prove that $\alpha$ is a periodic curve.

\begin{lemma} \label{rlopez:le3}
Under the hypothesis of this section and with the above notation, we have:
\begin{eqnarray*}
x(s+T)&=&x(s)+x(T)\\
z(s+T)&=&z(s)\\
\theta(s+T)&=&\theta(s)-2\pi
\end{eqnarray*}
\end{lemma}

\begin{proof}
This is a consequence of the uniqueness of solutions of
(\ref{rlopez:eq1})-(\ref{rlopez:eq2}). We only have to show that $z(T)=z_0$. But
this is a direct consequence of  the assumption (\ref{rlopez:c4}), that
$a^2+4b<0$ and  (\ref{rlopez:zeta}).
\end{proof}

As conclusion of Lemma  \ref{rlopez:le3},  we describe the
behavior of the coordinates functions of the profile curve $\alpha$
under the hypothesis (\ref{rlopez:c4}).  Due to
the monotonicity of $\theta$, let $T_1, T_2$ and $T_3$ be the points
in the interval $[0,T]$ such that the function $\theta$ takes the
values $-\pi/2, -\pi$ and $-3\pi/2$ respectively. In view of the
variation of the angle $\theta$ with the time coordinate $s$, it is
easy to verify the following  Table: \vspace*{4mm}

\begin{center}
\begin{tabular}{|c|c|c|c|}\hline
$s$ & $\theta$ & $x(s)$ & $z(s)$\\ \hline
$[0,T_1]$ & $[\frac{-\pi}{2},0]$ & \mbox{ increasing}
& \mbox{ decreasing}\\ \hline
$[T_1,T_2]$ & $[-\pi,\frac{-\pi}{2}]$ &\mbox{ decreasing} &
\mbox{ decreasing} \\ \hline
$[T_2,T_3]$ & $[\frac{-3\pi}{2},-\pi]$ &\mbox{ decreasing} &
\mbox{ increasing} \\ \hline
$[T_3,T]$ & $[-2\pi,\frac{-3\pi}{2}]$ &\mbox{ increasing} &
\mbox{ increasing}\\
\hline
\end{tabular}
\end{center}
\vspace*{4mm}

\begin{figure}[htbp]
\begin{center}\includegraphics[width=6cm]{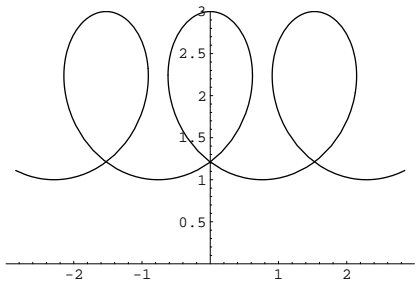}\end{center}
 \caption{The generating curve of a  rotational hyperbolic surfaces, with
 $a=-b=2$. Here  $z_0=3$. The curve $\alpha$ is periodic
with self-intersections. }\label{rlopez:fig3}
\end{figure}

\begin{theorem} \label{rlopez:t4}
Let $\alpha=\alpha(s)=(x(s),0,z(s))$ be the profile curve of a
rotational hyperbolic  surface $S$ in $\r^3$  where $\alpha$ is  the
solution of (\ref{rlopez:eq1})-(\ref{rlopez:eq2}). Assume that the initial
condition on $z_0$ satisfies $z_0>\frac{-2b}{a}$.  Then (see Fig. \ref{rlopez:fig3})
\begin{enumerate}
\item The curve $\alpha$ is invariant by the group of translations
in the $x$-direction given by the vector $(x(T),0,0)$.
\item In the period $[0,T]$ of $z$ given by Lemma \ref{rlopez:le3},  the function $z=z(s)$ presents one maximum at $s=0$
and one minimum at $s=T_2$.  Moreover, $\alpha$ is symmetric  with
respect to the vertical line at $x=0$ and $x=x(T_2)$.
\item The height function of  $\alpha$, that is, $z=z(s)$,  is periodic.
\item The curve $\alpha$ has self-intersections and its curvature has constant sign.
\item The part of $\alpha$ between the maximum and the minimum satisfies that the function $z(s)$ is
strictly decreasing with exactly one vertical point. Between this
minimum and the next maximum, $z=z(s)$ is strictly increasing  with
exactly one vertical point.
\item The velocity $\alpha'$  turns around the origin.
\end{enumerate}
\end{theorem}

\begin{theorem} \label{rlopez:t42} Let $S$ be a rotational hyperbolic  surface in $\r^3$ whose profile
curve $\alpha$ satisfies the hypothesis of  Theorem \ref{rlopez:t4}. Then
$S$ has the following properties:
\begin{enumerate}
\item The surface has self-intersections.
\item The surface is periodic with infinite vertical symmetries.
\item The surface is  complete.
\item The part of $\alpha$ between two consecutive vertical
points and containing a maximum corresponds with points of $S$ with
positive Gaussian curvature; on the other hand, if  this part
contains a minimum,  the Gaussian curvature  is negative.
\end{enumerate}
\end{theorem}

As it was announced in Theorem \ref{rlopez:t2}, and in order to distinguish from  the surfaces of negative constant Gaussian
curvature, we conclude from Theorem \ref{rlopez:t42} the following

\begin{corollary}\label{rlopez:co}There exists a one-parameter family of
rotational hyperbolic linear Weingarten surfaces that are complete
and with self-intersections in $\r^3$. Moreover, the generating
curves of these surfaces are periodic.
\end{corollary}

%%%%%%%%%%%%%%%%%%%%%%%%%%%%%%%%%%%%%%%%%
\section{Parabolic Weingarten surfaces in $\h^3$}\label{rlopez:sec4}
%%%%%%%%%%%%%%%%%%%%%%%%%%%%%%%%%%%%%%%%%%%%%%%%%%

A parabolic group of isometries of hyperbolic space $\h^3$ is
formed by isometries that leave fix one double point of the ideal
boundary $\s^2_{\infty}$  of $\h^3$.  A surface $S$ in $\h^3$ is called a
{\it parabolic surface}  if it is invariant by a group of
parabolic isometries. A  parabolic surface $S$ is determined by
a generating curve  $\alpha$ obtained  by the intersection of $S$
with any geodesic plane orthogonal to the orbits of the group.

We consider the upper half-space model of  $\h^3$, namely,
${\h}^{3}=:{\r}^3_+=\{(x,y,z)\in{\r}^{3};z>0\}$ equipped with the metric
$\langle,\rangle=\frac{dx^{2}+dy^{2}+dz^{2}}{z^{2}}$.
The ideal boundary $\s^2_{\infty}$ of $\h^3$ is identified with the one point compactification of the plane $\Pi\equiv\{z=0\}$, that is,
$\s^2_{\infty}=\Pi\cup\{\infty\}$. In what follows, we will  use the words  vertical or horizontal in the
usual affine sense of $\r^3_+$. Denote $L=\s^2_{\infty}\cap \{y=0\}$.

Let $G$ be a parabolic group of isometries. In the upper half-space model, we  take the point  $\infty\in \s^2_{\infty}$ as the point that fixes $G$.
 Then the group $G$ is defined by the horizontal (Euclidean) translations in the
direction of a  horizontal vector $\xi$ with $\xi\in\Pi$ which can be assumed $\xi=(0,1,0)$.

 A (parabolic) surface $S$ invariant by $G$ parametrizes as $X(s,t)=(x(s),t,z(s))$, where $t\in\r$
and the curve $\alpha=(x(s),0,z(s))$, $s\in I\subset\r$, is assumed to be parametrized by the Euclidean arc-length. The curve $\alpha$ is the generating curve of $S$. We write $\alpha'(s)=(\cos\theta(s),0,\sin\theta(s))$,
for a certain differentiable function $\theta$, where the derivative $\theta'(s)$ is
the Euclidean curvature of $\alpha$. With respect to the unit normal vector $N(s,t)=(-\sin\theta(s),0,\cos\theta(s))$, the principal curvatures are
$$\kappa_1(s,t)=z(s)\theta'(s)+\cos\theta(s),\hspace*{.5cm}\kappa_2(s,t)=\cos\theta(s).$$
The relation $aH+bK=c$ writes then
\begin{equation}\label{rlopez:w41}
\left(\frac{a}{2}+b\cos\theta(s)\right)z(s)\theta'(s)+a\cos\theta(s)-b\sin^2\theta(s)=c.
\end{equation}
We consider  initial conditions
\begin{equation}\label{rlopez:w42}
x(0)=0,\hspace*{.5cm}z(0)=z_0>0,\hspace*{.5cm}\theta(0)=0.
\end{equation}
Then any solution $\{x(s),z(s),\theta(s)\}$ satisfies properties of symmetry which are consequence of the uniqueness of solutions of an O.D.E. For example, the solution is symmetric with respect to the vertical straight line $x=0$. Using uniqueness again, we infer immediately

\begin{proposition}\label{rlopez:pr41}
Let $\alpha$ be  a solution of the initial value problem (\ref{rlopez:w41})-(\ref{rlopez:w42}) with $\theta(0)=\theta_0$.  If  $\theta'(s_0)=0$ at some  real number $s_0$, then $\alpha$ is parameterized by
$\alpha(s)=((\cos\theta_0) s,0,(\sin\theta_0) s+z_0)$, that is, $\alpha$ is a straight line  and the corresponding surface is a totally geodesic plane, an
equidistant surface or a horosphere.
\end{proposition}

In view of this proposition, we can assume that the function $\theta'(s)$ do not vanish, that is, $\theta$ is a monotonic function on $s$. At $s=0$, Equation (\ref{rlopez:w41}) is
 $$\theta'(0)=\frac{2}{z_0}\frac{c-a}{a+2b}.$$
This means that the study of solutions of (\ref{rlopez:w41})-(\ref{rlopez:w42}) must analyze a variety of cases depending on the sign of $\theta'(0)$. In this section, we are going to consider some cases in order to show techniques and some results. First, assume that $c\not=0$, which it can be assumed to be $c=1$.  Then we write (\ref{rlopez:w41}) as
\begin{equation}\label{rlopez:w233}
\theta'(s)=2\ \frac{1-a\cos\theta(s)+b\sin^2\theta(s)}{z(s)(a+2b\cos\theta(s))}.
\end{equation}
Our first result considers a case where it is possible to obtain explicit examples.

\begin{theorem}\label{rlopez:t41} Let $\alpha(s)=(x(s),0,z(s))$ be the generating curve of a
parabolic surface  $S$ in hyperbolic space $\h^3$ that satisfies
$aH+bK=1$ with initial conditions (\ref{rlopez:w42}). Assume   $a^2+4b^2+4b=0$. Then $\alpha$ describes an open of an
Euclidean circle in the $xz$-plane. \end{theorem}

\begin{proof}
Equation (\ref{rlopez:w41}) reduces into
$$-2bz(s)\theta'(s)=a+2b\cos\theta(s).$$
By differentiation with respect to $s$, we obtain
$z(s)\theta''(s)=0$, that is, $\theta'(s)$ is a constant function. Since $\theta'(s)$
describes the Euclidean curvature of $\alpha$, we conclude that
$\alpha$ parametrizes an Euclidean circle in the $xz$-plane and the assertion follows. This
circle may not to be completely included in the halfspace
${\r}^3_+$.
\end{proof}
 From now, we assume $a^2+4b^2+4b\not=0$.  Let us
denote by $(-\bar{s},\bar{s})$ the maximal domain of the solutions of (\ref{rlopez:w41})-(\ref{rlopez:w42}). By the monotonicity of
$\theta(s)$, let $\theta_1=\lim_{s\rightarrow\bar{s}}\theta(s)$.

\begin{figure}[htbp]\begin{center}
\includegraphics[width=4cm]{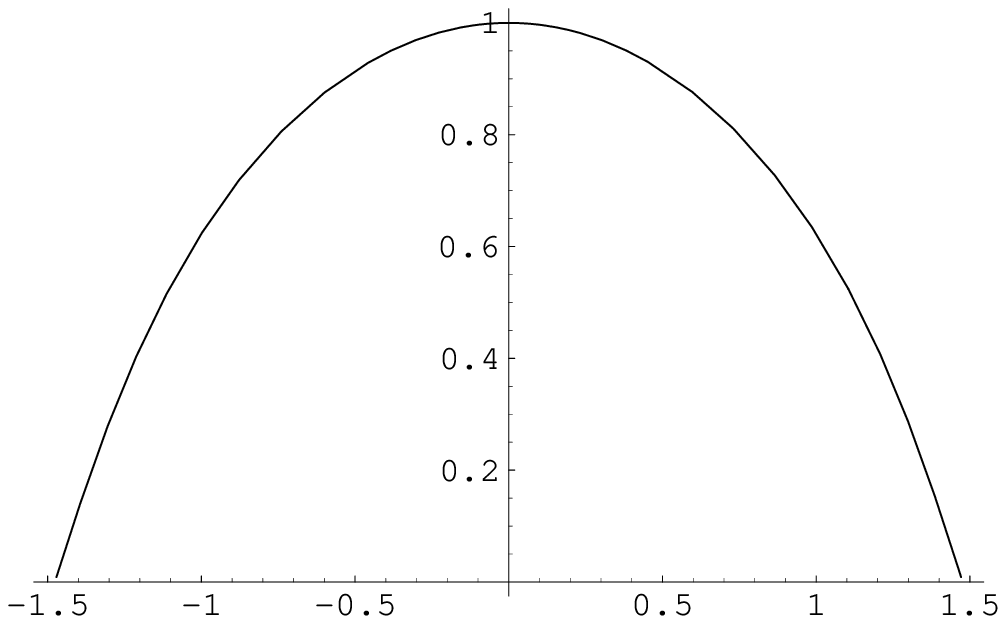}
\hspace*{2cm}\includegraphics[width=4cm]{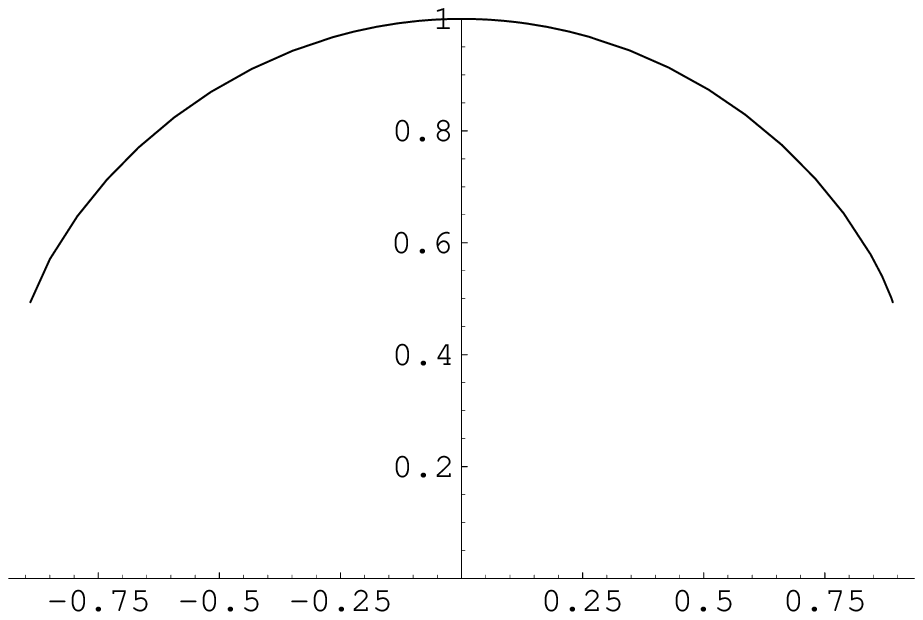}\end{center}
\begin{center}(a)\hspace*{5cm}(b)\end{center}
\caption{The generating curve of a parabolic surface with $aH+bK=1$,
with $0<a<1$ and $a+2b<0$.
 Here $z_0=1$ and $a=0.5$. In the case
(a), $b=-1$ and in the case (b), $b=-0.8$.}\label{rlopez:fig41}
\end{figure}

\begin{figure}[htbp]\begin{center}
\includegraphics[width=5cm]{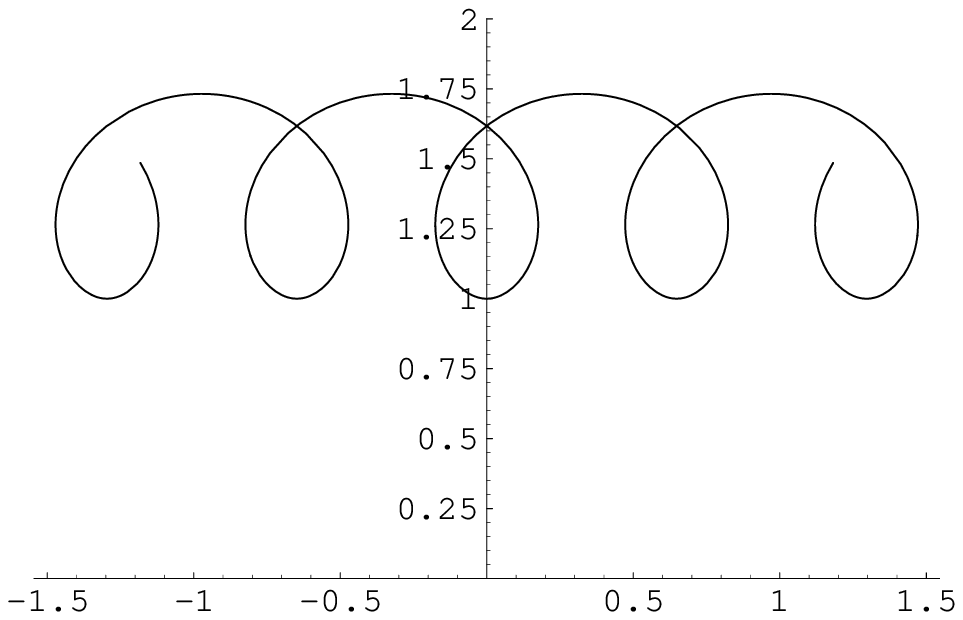}
\hspace*{2cm}\includegraphics[width=5cm]{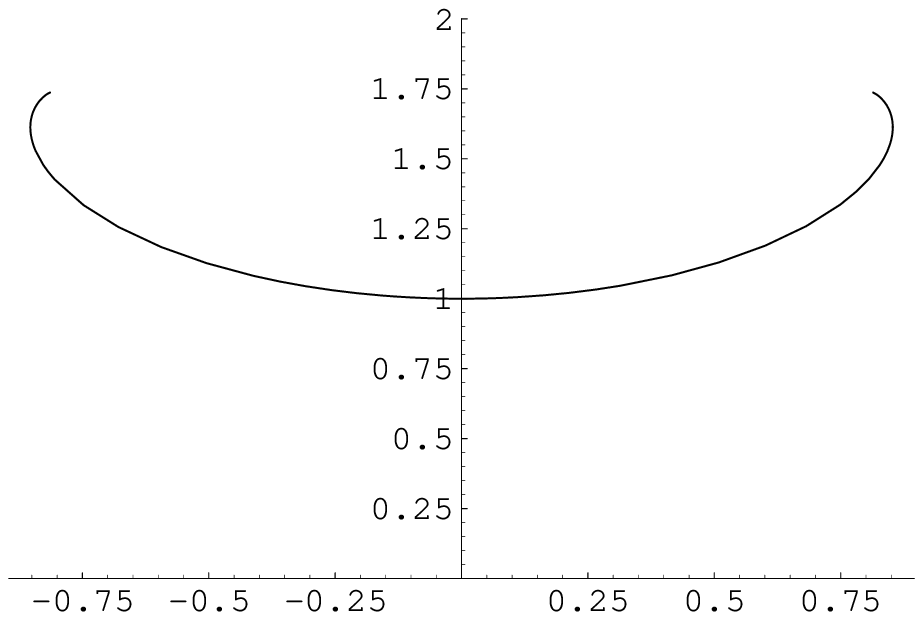}\end{center}
\begin{center}(a)\hspace*{8cm}(b)\end{center}
\caption{The generating curve of a parabolic surface with $aH+bK=1$,
with $0<a<1$ and $a+2b>0$.
 Here $z_0=1$ and $a=0.5$. In the case
(a), $b=-0.2$ and in the case (b), $b=0.3$.}\label{rlopez:fig42}
\end{figure}

\begin{theorem}[Case $0<a<1$]\label{rlopez:t43} Let $\alpha(s)=(x(s),0,z(s))$ be the generating curve of a
parabolic surface  $S$ in hyperbolic space $\h^3$ that satisfies $aH+bK=1$ with initial conditions (\ref{rlopez:w42}).
Assume $b\not=0$ and  $0<a<1$.
\begin{enumerate}
\item If $a+2b<0$, $\alpha$ has one maximum and $\alpha$ is a
concave (non-entire) vertical graph.  If $b<-(1+\sqrt{1-a^2})/2$, the surface $S$  is complete and intersects $\s^2_{\infty}$
 making an angle $\theta_1$ such that
$2\cos\theta_1-b\sin^2\theta_1=0$. The asymptotic boundary of $S$ is formed by two parallel straight lines. See Fig. \ref{rlopez:fig41} (a). If $-(1+\sqrt{1-a^2})/2<b<-a/2$,
then $S$ is not complete. See Fig. \ref{rlopez:fig41} (b).

\item Assume $a+2b>0$. If $a-2b>0$, then $S$ is complete and  invariant by a
group of translations in the $x$-direction. Moreover,  $\alpha$ has
self-intersections and it presents one maximum and one minimum in
each period.  See Fig. \ref{rlopez:fig42},  (a). If $a-2b\leq 0$, then $S$ is not complete. Moreover $\alpha$ is not a vertical graph  with  a minimum. See Fig.
\ref{rlopez:fig42}, (b).
\end{enumerate}
\end{theorem}


We point out that in each one of the cases of Theorem \ref{rlopez:t43}, we assert the existence of parabolic complete surfaces in $\h^3$ with the property $aH+bK=c$, such as it was announced in Theorem
\ref{rlopez:t3}.

\begin{proof}
The second derivative of $\theta''(s)$ satisfies
\begin{equation}\label{rlopez:dsegunda}
-\theta'(s)\sin\theta(s)\Big[b
\theta'(s)+\left(\frac{a}{2}+b\cos\theta(s)\right)\Big]+
\left(\frac{a}{2}+b\cos\theta(s)\right)z(s)\theta''(s)=0.
\end{equation}
\begin{enumerate}
\item Case  $a+2b<0$. Then $\theta'(0)<0$ and $\theta(s)$ is strictly decreasing.
If $\cos\theta(s)=0$ at some point $s$,
then (\ref{rlopez:w41}) gives $(a/2) z(s) \theta'(s)-b-1=0$. Thus, if
$b\geq -1$, $\cos\theta(s)\not=0$ and $-\pi/2<\theta(s)<\pi/2$.
In the case that $b<-1$ and as $a+2b\cos\theta(s)<0$,
 it follows from (\ref{rlopez:w41}) that  $a\cos\theta(s)-b\sin^2\theta(s)-1<0$ for any value of $s$, in particular,
 $\cos\theta(s)\not=0$. This proves that
$x'(s)=\cos\theta(s)\not=0$ and so,  $\alpha$ is a vertical graph on $L$. This
graph is concave since $z''(s)=\theta'(s)\cos\theta(s)<0$.
Moreover, this implies that $\bar{s}<\infty$ since on the contrary, and as $z(s)$ is
decreasing with $z(s)>0$, we would have $z'(s)\rightarrow 0$,
that is, $\theta(s)\rightarrow 0$: contradiction.

For $s>0$, $z'(s)=\sin\theta(s)<0$ and $z(s)$
 is strictly decreasing. Set $z(s)\rightarrow z(\bar{s})\geq 0$. The two roots of
$4b^2+4b+a^2=0$ on $b$ are $b=-\frac12(1\pm\sqrt{1-a^2})$. Moreover,
and from $a+2b<0$, we have
$$-\frac12(1+\sqrt{1-a^2})<\frac{-a}{2}<-\frac12(1-\sqrt{1-a^2}).$$
\begin{enumerate}
\item Subcase $b< -(1+\sqrt{1-a^2})/2$. Under this assumption,
$a^2+4b^2+4b>0$. Since $a<1$, we
obtain
\begin{equation} \label{rlopez:casouno}
a+2b\cos\theta(s)<-\sqrt{a^2+4b^2+4b}.
\end{equation}
If $z(\bar{s})>0$, then
$\lim_{s\rightarrow\bar{s}}\theta'(s)=-\infty$. In view of (\ref{rlopez:w233}) we have $a+2b \cos\theta(\bar{s})=0$: contradiction with
(\ref{rlopez:casouno}). Hence, $z(\bar{s})=0$ and $\alpha$ intersects $L$
with an angle
 $\theta_1$
satisfying $a\cos\theta_1-b\sin^2\theta_1-1=0$.
\item Subcase $-(1+\sqrt{1-a^2})/2< b<-a/2$. Now $a^2+4b^2+4b<0$. The function
 $1-a\cos\theta(s)+b\sin^2\theta(s)$ is strictly decreasing and its value
 at $\bar{s}$ satisfies $\cos\theta(s)>-a/2b$. Thus
 \begin{equation}\label{rlopez:casodos}
 1-a\cos\theta(s)+b\sin^2\theta(s)\geq \frac{a^2+4b^2+4b}{4b}>0.
 \end{equation}
  Assume $z(\bar{s})=0$.
Then (\ref{rlopez:w233}) and (\ref{rlopez:casodos}) imply that
$\theta'(\bar{s})=-\infty$. Combining
(\ref{rlopez:w233}) and (\ref{rlopez:dsegunda}), we have
 $$\frac{\theta''(s)}{\theta'(s)^2}
 =\frac{b\sin\theta(s)}{z(s)\left(\frac{a}{2}+b\cos\theta(s)\right)}+
\frac{\sin\theta(s)\left(\frac{a}{2}+b\cos\theta(s)\right)}{1-a\cos\theta(s)+b\sin^2\theta(s)}.
$$
From this expression and as  $\sin\theta(\bar{s})\not=0$, we conclude
$$\lim_{s\rightarrow\bar{s}}\frac{\theta''(s)}{\theta'(s)^2}=-\infty.$$
On the other hand, using L'H\^opital rule, we have
$$\lim_{s\rightarrow\bar{s}}z(s)\theta'(s)=
\lim_{s\rightarrow\bar{s}}-\frac{\sin\theta(s)}{\frac{\theta''(s)}{\theta'(s)^2}}=0.$$
By letting $s\rightarrow\bar{s}$ in (\ref{rlopez:w233}), we obtain a contradiction.
 Thus, $z(\bar{s})>0$. This means
that $\lim_{s\rightarrow\bar{s}}\theta'(s)=-\infty$ and it follows from
(\ref{rlopez:w41}) that
$$\lim_{s\rightarrow\bar{s}}\left(\frac{a}{2}+b\cos\theta(s)\right)=0.$$

\end{enumerate}
\item Case  $a+2b>0$. Then $\theta'(0)>0$ and $\theta(s)$ is strictly increasing.
We distinguish  two possibilities:
\begin{enumerate}
\item Subcase $a-2b>0$. We prove that $\theta(s)$ reaches the value $\pi$. On the
contrary, $\theta(s)<\pi$ and $z(s)$ is an increasing function. The
hypothesis $a-2b>0$ together $a+2b>0$ implies that
$a+2b\cos\theta(s)\geq \delta>0$ for some number $\delta$. From
(\ref{rlopez:w233}), $\theta'(s)$ is bounded and then $\bar{s}=\infty$. In
particular, $\lim_{s\rightarrow\infty}\theta'(s)=0$. As both $a-2b$
and $a+2b$ are positive numbers, the function $b\theta'(s)+(a+2b\cos
\theta(s))$ is positive near $\bar{s}=\infty$. Then using
(\ref{rlopez:dsegunda}),  $\theta''(s)$ is
positive for a certain value of $s$ big enough, which it is impossible. As conclusion, $\theta(s)$ reaches
the value $\pi$ at some $s=s_0$. By the symmetry properties of solutions of (\ref{rlopez:w41}),  $\alpha$ is
symmetric with respect to the line $x=x(s_0)$ and the velocity
vector of $\alpha$ rotates until to the initial position. This means
that $\alpha$ is invariant by a group of horizontal translations.
\item Subcase $a-2b\leq 0$. As $\theta'(s)>0$, Equation (\ref{rlopez:w233}) says that
$\cos\theta(s)\not=-1$, and so, $\theta(s)$ is bounded by
$-\pi<\theta(s)<\pi$.  As in the above subcase, if $\bar{s}=\infty$,
then $\theta'(s)\rightarrow 0$, and this is a contradiction. Then
$\bar{s}<\infty$ and $\lim_{s\rightarrow\bar{s}}\theta'(s)=\infty$.
Hence, $\cos\theta(\bar{s})=-a/(2b)$ and $\theta(s)$ reaches
the value $\pi/2$.
\end{enumerate}
\end{enumerate}
\end{proof}

Address: \\
Departamento de Geometr\'{\i}a y Topolog\'{\i}a\\
Universidad de Granada, Spain\\
e mail: {\tt rcamino@ugr.es}
\end{document}